\documentclass{article}[12pt]

\usepackage{latexsym}
\usepackage{amsmath}
\usepackage{amssymb}

\begin{document}

\title{SOME SPECIAL SOLUTIONS OF A NONLINEAR SYSTEM OF $4$ ORDINARY DIFFERENTIAL
EQUATIONS RECENTLY INTRODUCED TO INVESTIGATE THE EVOLUTION OF HUMAN RESPIRATORY VIRUS EPIDEMICS}

\author{Francesco Calogero$^{a,b}$\thanks{e-mail: francesco.calogero@roma1.infn.it},
\thanks{e-mail: francesco.calogero@uniroma1.it} , Andrea Giansanti $^{a,b}$\thanks{e-mail: andrea.giansanti@uniroma1.it}
 \thanks{e-mail: andrea.giansanti@roma1.infn.it}
 , Farrin Payandeh$^c$\thanks{e-mail: farrinpayandeh@yahoo.com}
 \thanks{e-mail: f$\_$payandeh@pnu.ac.ir}}

\maketitle   \centerline{\it $^{a}$Physics Department, University of
Rome "La Sapienza", Rome, Italy}

\maketitle   \centerline{\it $^{b}$INFN, Sezione di Roma 1}

\maketitle

\maketitle   \centerline{\it $^{c}$Department of Physics, Payame
Noor University, PO BOX 19395-3697 Tehran, Iran}

\maketitle

\begin{abstract}

A system of $4$ nonlinearly-coupled Ordinary Differential Equations has been
recently introduced to investigate the evolution of human respiratory virus
epidemics. In this paper we point out that some \textit{explicit }solutions
of that system can be obtained by \textit{algebraic} operations, provided
the parameters of the model satisfy certain \textit{constraints}.

\end{abstract}

\section{\textbf{Introduction}}

The following system of $4$ nonlinearly-coupled Ordinary Differential
Equations (ODEs) has been recently introduced to investigate the evolution
of human respiratory virus epidemics \cite{RWK2021}:
\begin{subequations}
\label{Pandemics}
\begin{equation}
\overset{\cdot }{\tilde{x}}_{1}=-k_{D}\tilde{x}_{1}+\alpha k_{R}\left(
\tilde{x}_{3}+\tilde{x}_{4}\right) ~,  \label{x1dotBis}
\end{equation}%
\begin{equation}
\overset{\cdot }{\tilde{x}}_{2}=k_{B}\tilde{x}_{1}+\left[ k_{B}-k_{D}-f%
\left( \tilde{x}_{1},\tilde{x}_{2},\tilde{x}_{3},\tilde{x}_{4}\right) \right]
\tilde{x}_{2}+\left[ k_{B}+\left( 1-\alpha \right) k_{R}\right] \left(
\tilde{x}_{3}+\tilde{x}_{4}\right) ~,  \label{x2dotBis}
\end{equation}%
\begin{equation}
\overset{\cdot }{\tilde{x}}_{3}=f\left( \tilde{x}_{1},\tilde{x}_{2},\tilde{x}%
_{3},\tilde{x}_{4}\right) \tilde{x}_{2}-\left( k_{R}+k_{D}+k_{P}\right)
\tilde{x}_{3}~,  \label{x3dotBis}
\end{equation}%
\begin{equation}
\overset{\cdot }{\tilde{x}}_{4}=k_{P}\tilde{x}_{3}-\left(
k_{R}+k_{D}+k_{DV}\right) \tilde{x}_{4}~,  \label{x4dotBis}
\end{equation}%
where%
\begin{equation}
f\left( \tilde{x}_{1},\tilde{x}_{2},\tilde{x}_{3},\tilde{x}_{4}\right) =%
\frac{k_{I}\left( \tilde{x}_{3}+\beta \tilde{x}_{4}\right) }{\tilde{x}_{1}+%
\tilde{x}_{2}+\tilde{x}_{3}+\beta \tilde{x}_{4}}~.  \label{fxxxxBis}
\end{equation}

\textbf{Notation}. We maintained the original notation of \cite{RWK2021},
except for the following replacement of the $4$ dependent variables $I\left(
t\right) $ (number of \textit{Immune} hosts),$~S\left( t\right) $ (number of
\textit{Susceptible} hosts),$~A\left( t\right) $ (number of \textit{%
Asymptomatic} and infectious hosts),$~C\left( t\right) $ (number of \textit{%
Symptomatic} and infectious hosts) used there, and the use of a superimposed
dot (instead of an appended prime) to denote differentiation with respect to
the dependent variable $t$ ("time"):
\end{subequations}
\begin{equation}
I\left( t\right) =\tilde{x}_{1}\left( t\right) ~,~~S\left( t\right) =\tilde{x%
}_{2}\left( t\right) ~,~~A\left( t\right) =\tilde{x}_{3}\left( t\right)
~,~~C\left( t\right) =\tilde{x}_{4}\left( t\right) ~;~~\overset{\cdot }{%
\tilde{x}}\left( t\right) \equiv d\tilde{x}(t)/dt~.
\end{equation}%
Note that in the following we occasionally \textit{omit} to indicate \textit{%
explicitly} the $t$-dependence of the dependent variables: and see below
\textbf{Remark 1.1 }for our (notational) motivation of the tilde
superimposed on these coordinates $\tilde{x}_{n}\left( t\right) $. $%
\blacksquare $

For the epidemiological significance of this model see \cite{RWK2021}, as
well as for references to analogous models.

\textbf{Remark 1.1}. The system of $4$ nonlinearly-coupled ODEs (\ref%
{Pandemics}) features the $8$ \textit{a priori arbitrary }(of course
time-independent) \textit{parameters }$k_{D},$ $k_{R},$ $k_{B},$ $k_{P},$ $%
k_{DV,}$ $k_{I},$ $\alpha ,$ $\beta ,$ for whose epidemiological
significance we refer to \cite{RWK2021}; in this paper we focus mainly on
some mathematical properties of this system, so we generally assume that
these are $8$ \textit{a priori arbitrary (}possibly even \textit{complex)
numbers, }although we shall comment occasionally on the relevance of such
mathematical treatment on the \textit{epidemic} model (when these numbers
are \textit{positive real} numbers).

One observation which is relevant for the mathematical discussion of this
model (\ref{Pandemics})---which we think is reasonable to state at the very
beginning of this paper---is to note that the parameter $k_{D}$ plays a
relatively trivial role in this system, because it can be altogether
eliminated from it via the following very simple change of dependent
variables:
\begin{subequations}
\label{xntildexn}
\begin{equation}
\tilde{x}_{n}\left( t\right) =x_{n}\left( t\right) ~\exp \left(
-k_{D}t\right) ~,~~~n=1,2,3,4~,  \label{xtildex}
\end{equation}%
implying of course%
\begin{equation}
\tilde{x}_{n}\left( 0\right) =x_{n}\left( 0\right) ~,~~~n=1,2,3,4~;
\end{equation}%
indeed the system of ODEs satisfied by the $4$ variables $x_{n}\left(
t\right) $ is then \textit{identical} with the original system (\ref%
{Pandemics}), except for the elimination of the parameter $k_{D}$:
\end{subequations}
\begin{subequations}
\label{eqs1}
\begin{equation}
\dot{x}_{1}=\alpha k_{R}\left( x_{3}+x_{4}\right) ~,  \label{x1dot}
\end{equation}%
\begin{equation}
\dot{x}_{2}=k_{B}x_{1}+\left[ k_{B}-f\left( x_{1},x_{2},x_{3},x_{4}\right) %
\right] x_{2}+\left[ k_{B}+\left( 1-\alpha \right) k_{R}\right] \left(
x_{3}+x_{4}\right) ~,  \label{x2dot}
\end{equation}%
\begin{equation}
\dot{x}_{3}=f\left( x_{1},x_{2},x_{3},x_{4}\right) x_{2}-\left(
k_{R}+k_{P}\right) x_{3}~,  \label{x3dot}
\end{equation}%
\begin{equation}
\dot{x}_{4}=k_{P}x_{3}-\left( k_{R}+k_{DV}\right) x_{4}~,  \label{x4dot}
\end{equation}%
where of course now%
\begin{equation}
f\left( x_{1},x_{2},x_{3},x_{4}\right) =\frac{k_{I}\left( x_{3}+\beta
x_{4}\right) }{x_{1}+x_{2}+x_{3}+\beta x_{4}}~.~~~\blacksquare  \label{fxxxx}
\end{equation}

Hence hereafter we shall mainly deal with this, marginally simpler, system (%
\ref{eqs1}).

In the following \textbf{Section 2} we investigate a very simple solution of
this model (\ref{eqs1}), characterized by the fact that the $4$ components $%
x_{n}\left( t\right) $ of this solution all evolve proportionally to the
same exponential function of time, $\exp \left( \mu t\right) $, with $\mu $
an appropriate parameter determined in terms of the parameters of the model;
implying that the quantity $f\left( x_{1},x_{2},x_{3},x_{4}\right) $ is
\textit{time-independent} (see (\ref{fxxxx})), hence that the system (\ref%
{eqs1})---for this class of solutions---reduces to a \textit{linear} system
of $4$ ODEs.

In the subsequent \textbf{Section 3} we discuss the somewhat less simple
solutions characterized by the presumably more interesting requirement that
each of the $4$ components $x_{n}\left( t\right) $ of the solution be
\textit{linear} combinations---with \textit{time-independent}
coefficients---of $2$ exponential functions of time, $\exp \left( \mu
_{1}t\right) $ and $\exp \left( \mu _{2}t\right) ,$ and moreover that the
quantity $f\left( x_{1},x_{2},x_{3},x_{4}\right) ,$ see (\ref{fxxxx}), be
again \textit{time-independent}. The related restrictions on the parameters
of the model and the initial-data of this solution are also \textit{%
explicitly} determined, up to \textit{algebraic} operations.

The subsequent \textbf{Section 4} outlines the \textit{analogous} treatments
when the solution being identified is the sum of $3$, or $4$, exponentials.

A final \textbf{Section 5} concludes the paper, by mentioning its
applicative relevance and possible further developments of the approach used
in this paper.

\bigskip

\section{A very simple solution}

The right-hand sides of the $4$ ODEs (\ref{eqs1}) are all \textit{homogeneous%
} of degree $1$ in the $4$ dependent variables $x_{n}\left( t\right) .$ This
implies a well-known (see for instance \cite{CP2021}) consequence, which can
be stated as the following

\textbf{Proposition 2.1}. The system of $4$ ODEs (\ref{eqs1}) features the
\textit{simple} \textit{explicit} solution
\end{subequations}
\begin{equation}
x_{n}\left( t\right) =x_{n}\left( 0\right) ~\exp \left( \mu t\right)
~,~~~n=1,2,3,4~,  \label{SolProp1}
\end{equation}%
where $x_{n}\left( 0\right) $ are clearly the $4$ \textit{initial values} of
the $4$ dependent variables $x_{n}\left( t\right) $ and $\mu $ is an \textit{%
a priori arbitrary} time-independent parameter, provided these $5$
quantities---i. e., $x_{n}\left( 0\right) $ and $\mu $, together with the $7$
parameters of the model (\ref{eqs1})---satisfy (as it were, \textit{a
posteriori}) the following $4$ \textit{algebraic} equations:
\begin{subequations}
\label{Exn}
\begin{equation}
\mu x_{1}\left( 0\right) =\alpha k_{R}\left[ x_{3}\left( 0\right)
+x_{4}\left( 0\right) \right] ~,  \label{Ex1}
\end{equation}%
\begin{equation}
\mu x_{2}\left( 0\right) =k_{B}x_{1}\left( 0\right) +\left[ k_{B}-f\left(
0\right) \right] x_{2}\left( 0\right) +\left[ k_{B}+\left( 1-\alpha \right)
k_{R}\right] \left[ x_{3}\left( 0\right) +x_{4}\left( 0\right) \right] ~,
\label{Ex2}
\end{equation}%
\begin{equation}
\mu x_{3}\left( 0\right) =f\left( 0\right) x_{2}\left( 0\right) -\left(
k_{R}+k_{P}\right) x_{3}\left( 0\right) ~,  \label{Ex3}
\end{equation}%
\begin{equation}
\mu x_{4}\left( 0\right) =k_{P}x_{3}\left( 0\right) -\left(
k_{R}+k_{DV}\right) x_{4}\left( 0\right) ~,  \label{Ex4}
\end{equation}%
where of course (see (\ref{fxxxx}) and (\ref{SolProp1}))%
\begin{equation}
f\left( x_{1},x_{2},x_{3},x_{4}\right) \equiv f\left( 0\right) =\frac{k_{I}%
\left[ x_{3}\left( 0\right) +\beta x_{4}\left( 0\right) \right] }{%
x_{1}\left( 0\right) +x_{2}\left( 0\right) +x_{3}\left( 0\right) +\beta
x_{4}\left( 0\right) }~.~\blacksquare  \label{f(0)}
\end{equation}

The validity of this \textbf{Proposition 2.1} can be easily verified by
inserting the solution (\ref{SolProp1}) in the system (\ref{eqs1}) and by
then taking advantage of the conditions (\ref{Exn}).

Somewhat less trivial is to ascertain which are the \textit{constraints} on
the $4$ initial data $x_{n}\left( 0\right) $ and on the parameter $\mu $%
---by solving the system of algebraic equations (\ref{Exn})---when we
consider the model (\ref{eqs1}) for an \textit{arbitrary} assignment of its $%
7$ parameters $k_{R},k_{B},k_{P},k_{DV},k_{I},\alpha ,\beta $. Remarkably,
as we show below, this turns out to be \textit{explicitly} doable by purely
\textit{algebraic} operations.

Since all the $5$ eqs. (\ref{Exn}) are invariant under a \textit{common}
rescaling of the $4$ initial data $x_{n}\left( 0\right) $, it is convenient
to assume that one of them---say $x_{4}\left( 0\right) $---can be \textit{%
arbitrarily assigned}, and to focus on the \textit{ratios} of the other $3$
to that one, hence on the $3$ quantities
\end{subequations}
\begin{equation}
r_{m}=x_{m}\left( 0\right) /x_{4}\left( 0\right) ~,~~~x_{m}\left( 0\right)
=r_{m}x_{4}\left( 0\right) ~,~~~m=1,2,3~;  \label{rm}
\end{equation}%
thereby replacing the $5$ eqs. (\ref{Exn}) with the following $5$ equations:
\begin{subequations}
\label{Ern}
\begin{equation}
\mu ~r_{1}=\alpha k_{R}\left( r_{3}+1\right) ~,  \label{Er1}
\end{equation}%
\begin{equation}
\mu ~r_{2}=k_{B}r_{1}+\left[ k_{B}-F\left( r_{1},r_{2,}r_{3}\right) \right]
r_{2}+\left[ k_{B}+\left( 1-\alpha \right) k_{R}\right] \left(
r_{3}+1\right) ~,  \label{Er2}
\end{equation}%
\begin{equation}
\mu ~r_{3}=F\left( r_{1},r_{2,}r_{3}\right) r_{2}-\left( k_{R}+k_{P}\right)
r_{3}~,  \label{Er3}
\end{equation}%
\begin{equation}
\mu =k_{P}r_{3}-\left( k_{R}+k_{DV}\right) ~,  \label{Er4}
\end{equation}%
where of course (above and hereafter)%
\begin{equation}
F\left( r_{1},r_{2,}r_{3}\right) =k_{I}\left( r_{3}+\beta \right) /\left(
r_{1}+r_{2}+r_{3}+\beta \right) ~.  \label{F}
\end{equation}

It is now convenient---in order to get rid of the \textit{nonlinear}
function $F\left( r_{1},r_{2,}r_{3}\right) $---to sum the $2$ eqs. (\ref{Er2}%
) and (\ref{Er3}), getting thereby
\end{subequations}
\begin{equation}
\mu ~\left( r_{2}+r_{3}\right) =k_{B}+\left( 1-\alpha \right)
k_{R}+k_{B}r_{1}+k_{B}r_{2}+\left( k_{B}-k_{P}-\alpha k_{R}\right) r_{3}~.
\label{Er5}
\end{equation}

The $3$ eqs. (\ref{Er1}), (\ref{Er4}) and (\ref{Er5}) constitute now a
system of \textit{3} \textit{linear algebraic} equations for the $3$
unknowns $r_{1},r_{2},r_{3}$, which can be easily solved. Indeed from (\ref%
{Er4}) we get
\begin{subequations}
\label{rn}
\begin{equation}
r_{3}=\left( \mu +k_{DV}+k_{R}\right) /k_{P}~;  \label{r3}
\end{equation}%
then from (\ref{Er1}) and (\ref{r3}) we get%
\begin{equation}
r_{1}=\alpha k_{R}(\mu +k_{DV}+k_{P}+k_{R})/\left( \mu k_{P}\right) ~;
\label{r1}
\end{equation}%
and then from (\ref{Er5}), (\ref{r3}) and (\ref{r1}) we get
\begin{eqnarray}
&&r_{2}=-\frac{\mu +k_{DV}}{k_{P}}-\frac{\mu -k_{B}+k_{DV}}{\mu -k_{B}}
\notag \\
&&-\frac{\left[ \left( 1+\alpha \right) \mu +\alpha \left(
k_{DV}+k_{P}\right) \right] k_{R}+\alpha \left( k_{R}\right) ^{2}}{\mu k_{P}}%
~.  \label{r2}
\end{eqnarray}%
Note that these are \textit{explicit} expressions of the $3$ parameters $%
r_{m}$ in terms of the $6$ parameters $k_{R},$ $k_{B},$ $k_{P},$ $k_{DV},$ $%
k_{I},$ $\alpha $ of the system (\ref{eqs1}), and moreover of the parameter $%
\mu $ featured by the solution (\ref{SolProp1}) (where of course now $%
x_{m}\left( 0\right) =r_{m}x_{4}\left( 0\right) $ for $m=1,2,3,$ with $%
x_{4}\left( 0\right) $ remaining as a \textit{free} parameter).

Our remaining task in order to get the special solution (\ref{SolProp1}) of
the system (\ref{eqs1}) is to ascertain the permitted values of the
parameter $\mu $, as implied by inserting the following expression of $%
F(r_{1},r_{2},r_{3})$ (obtained by inserting the $3$ expressions (\ref{rn})
of $r_{1},r_{2},r_{3}$ in (\ref{F})),
\end{subequations}
\begin{equation}
F\left( r_{1},r_{2}r_{3}\right) =\frac{-k_{I}\left( \mu -k_{B}\right) \left(
\mu +k_{DV}+\beta k_{P}+k_{R}\right) }{k_{P}\left[ \left( 1-\beta \right)
\left( \mu -k_{B}\right) +k_{DV}\right] }~,  \label{FFF}
\end{equation}%
into any one of the $2$ eqs. (\ref{Er2}) or (\ref{Er3}).\textbf{\ }This
yields the following \textit{algebraic} equation of degree $4$ (hence
\textit{explicitly solvable}) for the quantity $\mu $:
\begin{subequations}
\label{eq8c}
\begin{equation}
\sum_{k=0}^{4}\left( c_{k}~\mu ^{k}\right) =0~,  \label{eqE}
\end{equation}%
with the following definitions of the $5$ parameters $c_{k}$:%
\begin{equation}
c_{4}=k_{I}-k_{P}+\beta k_{P},  \label{ccc4}
\end{equation}%
\begin{eqnarray}
&&c_{3}=2k_{DV}k_{I}-2k_{DV}k_{P}+k_{I}k_{P}-\left( k_{P}\right)
^{2}+2k_{I}k_{R}-2k_{P}k_{R}  \notag \\
&&+\alpha k_{I}k_{R}-k_{B}\left[ k_{I}-(1-\beta )k_{P}\right] +\beta
k_{P}(k_{DV}+k_{I}+k_{P}+2k_{R})~,  \label{ccc3}
\end{eqnarray}%
\begin{eqnarray}
&&c_{2}=(k_{I}-k_{P})\left[ \left( k_{DV}\right) ^{2}+k_{R}(k_{P}+k_{R})%
\right] +\alpha k_{I}k_{R}(k_{P}+2k_{R})  \notag \\
&&+\beta k_{P}\left[ (k_{I}+k_{R})(k_{P}+k_{R})+\alpha k_{I}k_{R}\right]
\notag \\
&&+k_{DV}\left\{ \left( 2+\beta \right) k_{I}k_{P}+2k_{I}\left( k_{R}+\alpha
k_{R}\right) -k_{P}\left[ (2-\beta )k_{P}+(3-\beta )k_{R}\right] \right\}
\notag \\
&&+k_{B}\left\{ {}\right. (1-\beta
)k_{P}(k_{P}+2k_{R})+k_{DV}(-2k_{I}+k_{P}-\beta k_{P})  \notag \\
&&-k_{I}(k_{P}+2k_{R}+\alpha k_{R}+\beta k_{P})\left. {}\right\} ~,
\label{ccc2}
\end{eqnarray}%
\begin{eqnarray}
&&c_{1}=(k_{DV}+k_{R})\left[ k_{DV}k_{P}(k_{I}-k_{P}-k_{R})+\alpha
k_{I}k_{R}(k_{DV}+k_{P}+k_{R})\right]   \notag \\
&&+\beta k_{I}k_{P}\left[ k_{DV}k_{P}+\alpha k_{R}(k_{DV}+k_{P}+k_{R})\right]
\notag \\
&&+k_{B}\left\{ {}\right. -(k_{DV}+k_{R})\left[
k_{DV}k_{I}+(k_{I}-k_{P})(k_{P}+k_{R})\right]   \notag \\
&&-\alpha k_{I}k_{R}(2k_{DV}+k_{P}+2k_{R})  \notag \\
&&-\beta k_{P}\left[ (k_{I}+k_{R})(k_{P}+k_{R})+k_{DV}(k_{I}+k_{P}+k_{R})+%
\alpha k_{I}k_{R}\right] \left. {}\right\} ~,  \label{ccc1}
\end{eqnarray}

\begin{equation}
c_{0}=-\alpha k_{B}k_{I}k_{R}\left( k_{DV}+k_{P}+k_{R}\right) \left(
k_{DV}+k_{R}+\beta k_{P}\right) ~.  \label{ccc0}
\end{equation}

\textbf{Remark 2.2}. For completeness let us mention that the results
reported just above require the validity of the following \textit{%
inequalities}:
\end{subequations}
\begin{equation}
\mu \neq 0~,~~~k_{P}\neq 0~,~~~\left( 1-\beta \right) \left( k_{B}-\mu
\right) -k_{DV}\neq 0~.~~~\blacksquare
\end{equation}

Some of these expressions of the $5$ parameters $c_{k}$---see (\ref{eq8c}%
)---are rather cumbersome (albeit quite explicit), featuring the $7$ \textit{%
a priori arbitrary} parameters $k_{R},$ $k_{B},$ $k_{P},$ $k_{DV},$ $k_{I},$
$\alpha ,$ $\beta $ characterizing the system (\ref{eqs1}); and of course
much more complicated are the---in principle easily available---\textit{%
explicit} expressions of the $4$ roots $\mu _{n}$ ($n=1,2,3,4$) of the
fourth-degree equation (\ref{eqE}). We do not consider useful to report
these formulas in this paper; since analogous---much more
practical---formulas can be easily obtained from eq. (\ref{eqE}) in \textit{%
applicative} contexts, whenever the $7$ \textit{a priori arbitrary}
parameters $k_{R},$ $k_{B},$ $k_{P},$ $k_{DV},$ $k_{I},$ $\alpha ,$ $\beta $
have been assigned specific numerical values, entailing, via the explicit
expressions of the parameters $c_{k}$ written above (see (\ref{eq8c})), the
corresponding numerical values of these parameters $c_{k}$, to be then
inserted in (\ref{eqE}) before the standard task of solving this \textit{%
quartic} equation is performed.

\textbf{Remark 2.3}. Let us finally mention that clearly, by setting
\begin{equation}
\mu =k_{D}~,  \label{EeqkD}
\end{equation}%
one is looking---see (\ref{xtildex}) and (\ref{SolProp1})---at the \textit{%
equilibrium} solution
\begin{equation}
\tilde{x}_{n}\left( t\right) =\bar{x}_{n}~,~~~\overset{\cdot }{\bar{x}}%
_{n}\left( t\right) =0~,~~~n=1,2,3,4~,  \label{equi}
\end{equation}%
of the original pandemic system (\ref{Pandemics}), as given by the formulas
(implied via (\ref{rm}) by (\ref{rn}))
\begin{subequations}
\label{xnbar}
\begin{equation}
\bar{x}_{1}=\bar{x}_{4}\left[ \alpha
k_{R}(k_{D}+k_{DV}+k_{P}+k_{R})/(k_{D}k_{P})\right] ~,
\end{equation}%
\begin{eqnarray}
&&\bar{x}_{2}=\bar{x}_{4}\left\{ -\frac{k_{D}+k_{DV}}{k_{P}}-\frac{%
k_{D}-k_{B}+k_{DV}}{k_{D}-k_{B}}\right. \\
&&-\left. \frac{\left[ \left( 1+\alpha \right) k_{D}+\alpha \left(
k_{DV}+k_{P}\right) \right] k_{R}+\alpha \left( k_{R}\right) ^{2}}{k_{P}k_{D}%
}\right\} ~,
\end{eqnarray}%
\begin{equation}
\bar{x}_{3}=\bar{x}_{4}(k_{D}+k_{DV}+k_{R})/k_{P}~,
\end{equation}%
where $\bar{x}_{4}$ is of course an \textit{arbitrary} parameter. $%
\blacksquare $

\textbf{Remark 2.4}. Note that throughout this paper we assume that the $4$
roots $\mu _{n}$ of the \textit{quartic} algebraic equation (\ref{eqE}) are
all different among themselves. $\blacksquare $

To conclude this \textbf{Section 2}, let us mention that the \textit{special}
solutions (\ref{SolProp1}) are not very interesting in \textit{applicative}
contexts, since they imply that the $4$ dependent variables $x_{n}\left(
t\right) $ \textit{all} evolve in the \textit{same}, very \textit{simple},
manner. But fortunately---as shown below---it is also possible to identify
other---presumably more interesting---\textit{explicit} solutions of the
system of nonlinear ODEs (\ref{eqs1}).

\bigskip

\section{A less simple solution: the linear combination of $2$ exponentials}

In this \textbf{Section 3} we investigate the following class of solutions
of the system of ODEs (\ref{eqs1}):
\end{subequations}
\begin{subequations}
\begin{equation}
x_{n}\left( t\right) =a_{n1}\exp \left( \mu _{1}t\right) +a_{n2}\exp \left(
\mu _{2}t\right) ~,~~~n=1,2,3,4~,  \label{3xnt}
\end{equation}%
where $\mu _{1}$ and $\mu _{2}$ are $2$ \textit{different} roots of eq. (\ref%
{eq8c}); while corresponding values for the $8$ time-independent parameters $%
a_{n1}$ and $a_{n2}$ are obtained below.

\textbf{Remark 3.1}. Since there are $4$ (assumedly \textit{different})
solutions $\mu $ of the fourth-order algebraic eq. (\ref{eqE}), there are $%
4\cdot 3/2=6$ different assignments of the pair of values $\mu _{1},$ $\mu
_{2}$. Note moreover that, even if the $7$ parameters $k_{R},$ $k_{B},$ $%
k_{P},$ $k_{DV,}$ $k_{I},$ $\alpha ,$ $\beta $ of the system of ODEs (\ref%
{eqs1}) are \textit{all real} numbers (as is of course the case in the
pandemics case), the $4$ solutions $\mu _{n}$ of the fourth-order algebraic
eq. (\ref{eqE}) need \textit{not} be \textit{real} numbers; but if the $7$
parameters of the system of ODEs (\ref{eqs1}) are \textit{all real} numbers,
then \textit{non-real} solutions of the algebraic eq. (\ref{eqE}) must be
present in \textit{complex conjugate} pairs. $\blacksquare $

\textbf{Remark 3.2}. Note that we are again assuming, throughout this
\textbf{Section 3}, that the quantity $f\left(
x_{1},x_{2},x_{3},x_{4}\right) $ in the system (\ref{eqs1}) is \textit{%
time-independent}, hence equal to its value at the initial time $t=0$ (see (%
\ref{f(0)})); although this property, which was obvious in the treatment of
\textbf{Section 2 }(see (\ref{fxxxx}) and (\ref{SolProp1}))---and which is
essential to justify the existence of the subclass of solutions (\ref{3xnt}%
)---is now instead far from obvious: indeed conditions for it to
hold---involving the \textit{initial data} of these solutions, and also $1$
\textit{constraint} on the parameters of the system (\ref{eqs1})--- shall
have to be ascertained, see below. $\blacksquare $

The $4$ eqs. (\ref{3xnt}) involve of course the following $4$ relations
among the $8$ parameters $a_{n1}$ and $a_{n2}$ and the $4$ initial data $%
x_{n}\left( 0\right) $:%
\begin{equation}
x_{n}\left( 0\right) =a_{n1}+a_{n2}~,~~~n=1,2,3,4~.  \label{xn0abn}
\end{equation}

The assumption (see \textbf{Remark 3.2}) that the function $f\left(
x_{1},x_{2},x_{3},x_{4}\right) $ be \textit{time-independent} implies that
the system (\ref{eqs1}) is again a \textit{linear} system of $4$ ODEs with
\textit{time-independent} parameters; hence each of the $2$ exponential
functions in the right-hand side of the \textit{ansatz} (\ref{3xnt}) must
satisfy (as it were, \textit{separately}) the system (\ref{eqs1}). Therefore
each of the $2$ sets of $4$ parameters $a_{n1}$ and $a_{n2}$ must satisfy
the same requirements---see the $4$ eqs. (\ref{Exn})---satisfied by the
initial data $x_{n}\left( 0\right) $ in the treatment of the previous
\textbf{Section 2}; namely there must now hold the $8$\textbf{\ }relations
\end{subequations}
\begin{subequations}
\label{anel}
\begin{equation}
\mu _{\ell }a_{1\ell }=\alpha k_{R}\left( a_{3\ell }+a_{4\ell }\right)
~,~~~\ell =1,2~,
\end{equation}%
\begin{equation}
\mu _{\ell }a_{2\ell }=k_{B}a_{1\ell }+\left[ k_{B}-f\left( 0\right) \right]
a_{2\ell }+\left[ k_{B}+\left( 1-\alpha \right) k_{R}\right] \left[ a_{3\ell
}+a_{4\ell }\right] ~,~~~\ell =1,2~,
\end{equation}%
\begin{equation}
\mu _{\ell }a_{3\ell }=f\left( 0\right) a_{2\ell }-\left( k_{R}+k_{P}\right)
a_{3\ell }~,~~~\ell =1,2~,
\end{equation}%
\begin{equation}
\mu _{\ell }a_{4\ell }=k_{P}a_{3\ell }-\left( k_{R}+k_{DV}\right) a_{4\ell
}~,~~~\ell =1,2~;
\end{equation}%
where of course we again set
\begin{equation}
f\left( x_{1},x_{2},x_{3},x_{4}\right) \equiv f\left( 0\right) =\frac{k_{I}%
\left[ x_{3}\left( 0\right) +\beta x_{4}\left( 0\right) \right] }{%
x_{1}\left( 0\right) +x_{2}\left( 0\right) +x_{3}\left( 0\right) +\beta
x_{4}\left( 0\right) }~;
\end{equation}%
but now with the $4$ initial data $x_{n}\left( 0\right) $ related to the $2$
parameters $a_{n1}$ and $a_{n2}$ by the eq. (\ref{xn0abn}).

We can therefore now proceed in close analogy to the treatment of the
previous \textbf{Section 2}, introducing $6$ parameters $b_{m\ell }$ via the
following position:
\end{subequations}
\begin{equation}
a_{m\ell }=b_{m\ell }a_{4\ell }~,~~~b_{m\ell }=a_{m\ell }/a_{4\ell
}~,~~~m=1,2,3,~~\ell =1,2~.  \label{3abmel}
\end{equation}%
These $6$ parameters $b_{m\ell }$ ($m=1,2,3,~\ell =1,2$) are then \textit{%
explicitly} expressed in terms of the parameters of the system (\ref{eqs1})
as follows (see (\ref{rn})):
\begin{subequations}
\label{3bel}
\begin{equation}
b_{1\ell }=\alpha k_{R}(\mu _{\ell }+k_{DV}+k_{P}+k_{R})/\left( k_{P}\mu
_{\ell }\right) ~,~~~\ell =1,2~,
\end{equation}%
\begin{eqnarray}
&&b_{2\ell }=-\frac{\mu _{\ell }+k_{DV}}{k_{p}}-\frac{\mu _{\ell
}-k_{B}+k_{DV}}{\mu _{\ell }-k_{B}}  \notag \\
&&-\frac{\left[ \left( 1+\alpha \right) \mu _{\ell }+\alpha \left(
k_{DV}+k_{P}\right) \right] k_{R}+\alpha \left( k_{R}\right) ^{2}}{k_{P}\mu
_{\ell }}~,~~~\ell =1,2~,
\end{eqnarray}%
\begin{equation}
b_{3\ell }=(\mu _{\ell }+k_{DV}+k_{R})/k_{P}~,~~~\ell =1,2~;
\end{equation}%
with $a_{41}$ and $a_{42}$ remaining as $2$ \textit{free} parameters.

To complete the treatment of this case, it is necessary to identify the
\textit{constraints} on the parameters of the system (\ref{eqs1}) and on the
parameters of the solution under present consideration, see (\ref{3xnt}),
which are necessary in order to comply with the requirement---essential for
our treatment---that the quantity $h\left( t\right) $, related to the
quantity $f\left( x_{1},x_{2},x_{3},x_{4}\right) $, see (\ref{fxxxx}), by
the simple relation
\end{subequations}
\begin{subequations}
\begin{equation}
f\left( x_{1},x_{2},x_{3},x_{4}\right) =k_{I}h\left( t\right)
\end{equation}%
---hence reading as follows,%
\begin{equation}
h\left( t\right) =\frac{x_{3}\left( t\right) +\beta x_{4}\left( t\right) }{%
x_{1}\left( t\right) +x_{2}\left( t\right) +x_{3}\left( t\right) +\beta
x_{4}\left( t\right) }  \label{3ht}
\end{equation}%
---be \textit{time-independent}: therefore given in terms of the initial
data as follows:%
\begin{equation}
h\left( t\right) =h\left( 0\right) =\frac{x_{3}\left( 0\right) +\beta
x_{4}\left( 0\right) }{x_{1}\left( 0\right) +x_{2}\left( 0\right)
+x_{3}\left( 0\right) +\beta x_{4}\left( 0\right) }~.
\end{equation}

To fulfill this task, we now note that the solutions $x_{n}\left( t\right) $
under consideration in this \textbf{Section 3} are defined by the relations (%
\ref{3xnt}), hence their insertion in the definition (\ref{3ht}) implies the
following expression of $h\left( t\right) $:%
\begin{equation}
h\left( t\right) =\frac{a_{31}+\beta a_{41}+\left( a_{32}+\beta
a_{42}\right) \exp \left[ \left( \mu _{2}-\mu _{1}\right) t\right] }{%
a_{11}+a_{21}+a_{31}+\beta a_{41}+\left( a_{12}+a_{22}+a_{32}+\beta
a_{42}\right) \exp \left[ \left( \mu _{2}-\mu _{1}\right) t\right] }~.
\label{3hht}
\end{equation}

It is therefore easily seen that the requirement that this expression of $%
h\left( t\right) $ be \textit{time-independent} implies that the $8$
parameters $a_{n\ell }$ ($n=1,2,3,4$; $\ell =1,2$) satisfy the following
\textit{single constraint} on the\textit{\ }$8$ parameters $a_{n\ell }$:
\end{subequations}
\begin{equation}
\left( a_{12}+a_{22}\right) \left( a_{31}+\beta a_{41}\right) -\left(
a_{11}+a_{21}\right) \left( a_{32}+\beta a_{42}\right) =0~;
\label{3cConamel}
\end{equation}%
entailing then that%
\begin{equation}
h\left( t\right) =h\left( 0\right) =\frac{a_{31}+\beta a_{41}}{%
a_{11}+a_{21}+a_{31}+\beta a_{41}}=\frac{a_{32}+\beta a_{42}}{%
a_{12}+a_{22}+a_{32}+\beta a_{42}}~.
\end{equation}

By inserting the formulas (\ref{3abmel}) in (\ref{3cConamel}) we then get,
for the $6$ parameters $b_{m\ell }$ ($m=1,2,3;$ $\ell =1,2$), the following
\textit{single constraint}:
\begin{equation}
\left( b_{12}+b_{22}\right) \left( b_{31}+\beta \right) -\left(
b_{11}+b_{21}\right) \left( b_{32}+\beta \right) =0~;
\end{equation}%
and, via (\ref{3bel}), we finally get the following \textit{single
constraint }on the $7$ parameters of the system (\ref{eqs1}) for the
existence of the solution (\ref{3xnt}):
\begin{subequations}
\begin{eqnarray}
\left( 1-\beta \right) \left( k_{B}\right) ^{2}-k_{B}k_{DV}+k_{DV}\left(
k_{DV}+\beta k_{P}+k_{R}+\mu _{1}\right) &&  \notag \\
+\left[ k_{DV}+\left( 1-\beta \right) \mu _{1}\right] \mu _{2}-\left(
1-\beta \right) k_{B}\left( \mu _{1}+\mu _{2}\right) =0~, &&  \label{3Constr}
\end{eqnarray}%
provided there hold the following inequalities:%
\begin{equation}
k_{P}\neq 0~,~~~k_{B}\neq \mu _{1}~,~~\ k_{B}\neq \mu _{2}~,~\ ~\mu _{1}\neq
\mu _{2}~.
\end{equation}

\textbf{Remark 3.3}. Let us recall that there are in general $6$ different
versions of the \textit{constraint} (\ref{3Constr}) due to the $6$ different
possible selections of the $2$ roots $\mu _{1}$ and $\mu _{2}$ (see \textbf{%
Remark 3.1}); and that the simplicity of this formula (\ref{3Constr}) as
providing a \textit{constraint} on the $7$ parameters of the system (\ref%
{eqs1}) is somewhat misleading, due to the (\textit{explicit} but) \textit{%
quite complicated} dependence on these parameters of the solutions $\mu _{1}$
and $\mu _{2}$ of the \textit{fourth-degree} algebraic equation (\ref{eqE}).
However---as already mentioned above---all these complicated relations
(including those yielding the \textit{initial data} of the class of
solutions considered in this \textbf{Section 3}) become much more easily
managed whenever any $6$ of the $7$ \textit{a priori arbitrary} parameters
featured by the system (\ref{eqs1}) are assigned \textit{specific numerical}
values, so that the remaining task left is to ascertain the values of the $7$%
-th parameter implied by the constraint (\ref{3Constr}) (as well as those
characterizing the \textit{initial data} $x_{n}\left( 0\right) $ of the
class of solutions considered in this \textbf{Section 3}), thereby
identifying the corresponding class of systems (\ref{eqs1}) featuring the
simple \textit{explicit} solutions of type (\ref{3xnt}). $\blacksquare $

\bigskip

\section{Solutions which are the linear superposition of 3 or 4 exponentials}

In this \textbf{Section 4} we treat the subclass of solutions of the system (%
\ref{eqs1}) whose time-evolution is a \textit{linear superposition} of $3$
or $4$ exponentials.

\textbf{Remark 4.1}. Hence, throughout this \textbf{Section 4, }we assume
that the quantity $f\left( x_{1},x_{2},x_{3},x_{4}\right) $ in (\ref{eqs1})
is \textit{time-independent}; although this property is---as in \textbf{%
Section 3}: see for instance \textbf{Remark 3.2}---far from obvious: indeed
conditions for it to hold---involving the \textit{initial data} of these
subclass of solutions, as well as \textit{constraints} on the parameters of
the system (\ref{eqs1})---shall have to be ascertained, see below. $%
\blacksquare $

\bigskip

\subsection{Solutions which are the linear superposition of 3 exponentials}

In this \textbf{Subsection 4.1 }we investigate the following class of
solutions of the system of ODEs (\ref{eqs1}):
\end{subequations}
\begin{subequations}
\label{41xnt0}
\begin{equation}
x_{n}\left( t\right) =a_{n1}\exp \left( \mu _{1}t\right) +a_{n2}\exp \left(
\mu _{2}t\right) +a_{n3}\exp \left( \mu _{3}t\right) ~,~~~n=1,2,3,4~,
\label{41xnt}
\end{equation}%
where $\mu _{1},$ $\mu _{2},$ $\mu _{3}$ are $3$ different roots of the $4$%
th-degree algebraic eq. (\ref{eqE}); while corresponding values for the $12$
time-independent parameters $a_{n1,}$ $a_{n2},$ $a_{n3}$ are obtained below.

Of course these formulas (\ref{41xnt}) imply the following relations among
the $4$ initial data $x_{n}\left( 0\right) $ and the $12$ parameters $a_{nj}$
($n=1,2,3,4;$ $j=1,2,3$):
\begin{equation}
x_{n}\left( 0\right) =a_{n1}+a_{n2}+a_{n3}~,~~~n=1,2,3,4~.  \label{41xn0}
\end{equation}

\textbf{Remark 4.1.1}. Clearly symbols such as $x_{n},$ $\mu _{n\ell },$ $%
a_{n\ell }$ need not have the same significance nor the same values when
appearing in different \textbf{Sections} or \textbf{Subsections} of this
paper. But of course the statements made in \textbf{Remark 3.1} concerning
the possibility that not all the $4$ roots of the fourth-order algebraic eq.
(\ref{eqE}) be \textit{real} numbers are generally valid. $\blacksquare $

\textbf{Remark 4.1.2}. Since there are $4$ (assumedly \textit{different})
solutions $\mu $ of the fourth-order algebraic eq. (\ref{eqE}), there are $4$
different selections---from the \textit{quartet} of solutions $\mu _{n}$ of
eq. (\ref{eqE})---of the \textit{trio} of values $\mu _{1},$ $\mu _{2},$ $%
\mu _{3}$ in the \textit{ansatz} (\ref{41xnt}). $\blacksquare $

Let us now proceed in close analogy to the treatment provided in \textbf{%
Section 3}. Again we assume that the quantity $f$ in the right-hand sides of
the ODEs (\ref{x2dot}) and (\ref{x3dot}) is a \textit{time-independent }%
parameter, up to identifying below conditions on the parameters of the
system (\ref{eqs1}) and on the initial data of the solution (\ref{41xnt0})
under consideration which are \textit{sufficient} to guarantee---as it were,
\textit{a posteriori}---that this be the case; hence that the system of ODEs
(\ref{eqs1}) be equivalent to a system of $4$ \textit{linear} ODEs,
featuring independent solutions $a\exp (\mu t)$ each depending exponentially
on the independent variable $t$ (which can therefore be added without
loosing the property to satisfy the system of ODEs (\ref{eqs1})).

We\ thus obtain---in analogy to the $8$ relations (\ref{anel})---the
following $12$ relations:
\end{subequations}
\begin{subequations}
\label{41anel}
\begin{equation}
\mu _{j}a_{1j}=\alpha k_{R}\left( a_{3j}+a_{4j}\right) ~,~~~j=1,2,3~,
\end{equation}%
\begin{equation}
\mu _{j}a_{2j}=k_{B}a_{1j}+\left[ k_{B}-f\left( 0\right) \right] a_{2j}+%
\left[ k_{B}+\left( 1-\alpha \right) k_{R}\right] \left(
a_{3j}+a_{4j}\right) ~,~~~j=1,2,3~,
\end{equation}%
\begin{equation}
\mu _{j}a_{3j}=f\left( 0\right) a_{2j}-\left( k_{R}+k_{P}\right)
a_{3j}~,~~~j=1,2,3~,
\end{equation}%
\begin{equation}
\mu _{j}a_{4j}=k_{P}a_{3j}-\left( k_{R}+k_{DV}\right) a_{4j}~,~~~j=1,2,3~.
\end{equation}

Next we set (in analogy to (\ref{3abmel}))
\end{subequations}
\begin{equation}
a_{mj}=b_{mj}a_{4j}~,~~~b_{mj}=a_{mj}/a_{4j}~,~~~m=1,2,3,~~j=1,2,3~,
\label{41mj}
\end{equation}%
getting thereby the following $9$ relations:
\begin{subequations}
\label{41bel}
\begin{equation}
b_{1j}=\alpha k_{R}(\mu _{j}+k_{DV}+k_{P}+k_{R})/\left( k_{P}\mu _{j}\right)
~,~~~j=1,2,3~,
\end{equation}%
\begin{eqnarray}
&&b_{2j}=-\frac{\mu _{j}+k_{DV}}{k_{p}}-\frac{\mu _{j}-k_{B}+k_{DV}}{\mu
_{j}-k_{B}}  \notag \\
&&-\frac{\left[ \left( 1+\alpha \right) \mu _{j}+\alpha \left(
k_{DV}+k_{P}\right) \right] k_{R}+\alpha \left( k_{R}\right) ^{2}}{k_{P}\mu
_{j}}~,~~~j=1,2,3~,
\end{eqnarray}%
\begin{equation}
b_{3j}=(\mu _{j}+k_{DV}+k_{R})/k_{P}~,~~~j=1,2,3~;
\end{equation}%
with $a_{41}$, $a_{42}$, $a_{43}$ remaining as $3$ \textit{free} parameters.

We must now investigate the restrictions on the parameters $a_{mj}$ implied
by the requirement that the quantity $h\left( t\right) $ be \textit{%
time-independent}. The analogous formula to (\ref{3hht}) now reads as
follows:
\end{subequations}
\begin{subequations}
\begin{equation}
h\left( t\right) =numh\left( t\right) /denh\left( t\right) ~,
\end{equation}%
\begin{eqnarray}
&&numh\left( t\right) =a_{31}+\beta a_{41}+\left( a_{32}+\beta a_{42}\right)
\exp \left[ \left( \mu _{2}-\mu _{1}\right) t\right]  \notag \\
&&+\left( a_{33}+\beta a_{43}\right) \exp \left[ \left( \mu _{3}-\mu
_{1}\right) t\right] ~,
\end{eqnarray}%
\begin{eqnarray}
&&denh\left( t\right) =a_{11}+a_{21}+a_{31}+\beta a_{41}  \notag \\
&&+\left( a_{12}+a_{22}+a_{32}+\beta a_{42}\right) \exp \left[ \left( \mu
_{2}-\mu _{1}\right) t\right]  \notag \\
&&+\left( a_{13}+a_{23}+a_{33}+\beta a_{43}\right) \exp \left[ \left( \mu
_{3}-\mu _{1}\right) t\right] ~.
\end{eqnarray}

It is then easily seen that the condition (\ref{3cConamel}) is now replaced
by the following $2$ restrictions:
\end{subequations}
\begin{equation}
\left( a_{11}+a_{21}\right) \left( a_{3k}+\beta a_{4k}\right) =\left(
a_{1k}+a_{2k}\right) \left( a_{31}+\beta a_{41}\right) ~,~~~k=2,3~.
\end{equation}

Hence, after the change of parameters (\ref{41mj}), we get the following $2$
\textit{constraints,}%
\begin{equation}
\left( b_{11}+b_{21}\right) \left( b_{3k}+\beta \right) =\left(
b_{1k}+b_{2k}\right) \left( b_{31}+\beta \right) ~,~~~k=2,3~,
\end{equation}%
on the $9$ parameters $b_{mj}$ ($m=1,2,3;~j=1,2,3$). Which, via the
expressions (\ref{41bel}) of these $9$ parameters, entail the following $2$
\textit{constraints} on the original $7$ parameters of the system (\ref{eqs1}%
):
\begin{subequations}
\begin{eqnarray}
\left( 1-\beta \right) \left( k_{B}\right) ^{2}+\left( k_{DV}\right)
^{2}+\left( 1-\beta \right) \mu _{1}\mu _{2}+k_{DV}\left( \beta
k_{P}+k_{R}+\mu _{1}+\mu _{2}\right) &&  \notag \\
-k_{B}\left[ k_{DV}+\left( 1-\beta \right) \left( \mu _{1}+\mu _{2}\right) %
\right] =0~, &&
\end{eqnarray}
\begin{eqnarray}
k_{DV}\left( k_{DV}+\beta k_{P}+k_{R}+\mu _{1}\right) +\left[ k_{DV}+\left(
1-\beta \right) \mu _{1}\right] \mu _{3} &&  \notag \\
-\left( 1-\beta \right) k_{B}\left( \mu _{1}+\mu _{3}\right) +k_{B}\left[
\left( 1-\beta \right) k_{B}-k_{DV}\right] =0. &&
\end{eqnarray}%
~

While of course the initial data $x_{n}\left( 0\right) $ of the solution
under consideration in this \textbf{Subsection 4.1} are explicitly given by
the formulas (\ref{41xn0}) with (\ref{41mj}) and (\ref{41bel}).

\bigskip

\subsection{Solutions which are the linear superposition of $4$ exponentials}

The treatment in this \textbf{Subsection 4.2} is quite terse, since it is
quite analogous to that provided above in \textbf{Subsection 4.1}; hence we
only report the key formulas which play an analogous role to the key
formulas in \textbf{Subsection 4.1}.

Instead of (\ref{41xnt0}) we now have
\end{subequations}
\begin{subequations}
\label{42xnt0}
\begin{equation}
x_{n}\left( t\right) =\sum_{q=1}^{4}\left[ a_{nq}\exp \left( \mu
_{q}t\right) \right] ~,~~~n=1,2,3,4~,  \label{42xnt}
\end{equation}%
\begin{equation}
x_{n}\left( 0\right) =\sum_{q=1}^{4}\left( a_{nq}\right) ~,~~~n=1,2,3,4~.
\label{42xn0}
\end{equation}

In place of the $12$ equations (\ref{41anel}) we now have the following $16$
equations:
\end{subequations}
\begin{subequations}
\label{42anq}
\begin{equation}
\mu _{q}a_{1q}=\alpha k_{R}\left( a_{3q}+a_{4q}\right) ~,~~~q=1,2,3,4~,
\label{42a1q}
\end{equation}%
\begin{eqnarray}
\mu _{q}a_{2q} &=&k_{B}a_{1q}+\left[ k_{B}-f\left( 0\right) \right] a_{2q}+%
\left[ k_{B}+\left( 1-\alpha \right) k_{R}\right] \left(
a_{3q}+a_{4q}\right) ~,  \notag \\
q &=&1,2,3,4~,  \label{42a2q}
\end{eqnarray}%
\begin{equation}
\mu _{q}a_{3q}=f\left( 0\right) a_{2q}-\left( k_{R}+k_{P}\right)
a_{3q}~,~~~q=1,2,3,4~,  \label{42a3nq}
\end{equation}%
\begin{equation}
\mu _{q}a_{4q}=k_{P}a_{3q}-\left( k_{R}+k_{DV}\right) a_{4q}~,~~~q=1,2,3,4~.
\label{42a4q}
\end{equation}

Likewise, in place of the $9$ equations (\ref{41mj}), we now write the $12$
relations
\end{subequations}
\begin{equation}
a_{jq}=b_{jq}a_{4j}~,~~~b_{jq}=a_{jq}/a_{4j}~,~~~j=1,2,3~,~~~q=1,2,3,4~,
\label{42abmq}
\end{equation}%
getting thereby, from (\ref{42anq}), the following $12$ relations (analogous
to (\ref{41bel})):
\begin{subequations}
\label{42bel}
\begin{equation}
b_{1q}=\alpha k_{R}(\mu _{q}+k_{DV}+k_{P}+k_{R})/\left( k_{P}\mu _{q}\right)
~,~~~q=1,2,3,4~,
\end{equation}%
\begin{eqnarray}
&&b_{2q}=-\frac{\mu _{q}+k_{DV}}{k_{p}}-\frac{\mu _{q}-k_{B}+k_{DV}}{\mu
_{q}-k_{B}}  \notag \\
&&-\frac{\left[ \left( 1+\alpha \right) \mu _{q}+\alpha \left(
k_{DV}+k_{P}\right) \right] k_{R}+\alpha \left( k_{R}\right) ^{2}}{k_{P}\mu
_{q}}~,~~~q=1,2,3,4~,
\end{eqnarray}%
\begin{equation}
b_{3q}=(\mu _{q}+k_{DV}+k_{R})/k_{P}~,~~~q=1,2,3,4~;
\end{equation}%
with $a_{41}$, $a_{42}$, $a_{43}$, $a_{44}$ remaining as $4$ \textit{free}
parameters.

Next come the \textit{constraints} on the parameters of the system (\ref%
{eqs1}) needed in order that the more general solution (\ref{42xnt}), when
inserted in the definition (\ref{fxxxx}) of the function $f\left( t\right) $%
, hence now reading
\end{subequations}
\begin{subequations}
\begin{equation}
f\left( t\right) =\frac{k_{I}\sum_{q=1}^{4}\left[ \left( a_{3q}+\beta
a_{4q}\right) \exp \left( \mu _{q}t\right) \right] }{\sum_{q=1}^{4}\left[
\left( a_{1q}+a_{2q}+a_{3q}+\beta a_{4q}\right) \exp \left( \mu _{q}t\right) %
\right] }~,
\end{equation}%
---or, equivalently, see (\ref{42abmq})---%
\begin{equation}
f\left( t\right) =\frac{k_{I}\sum_{q=1}^{4}\left[ \left( b_{3q}+1\right)
\exp \left( \mu _{q}t\right) \right] }{\sum_{q=1}^{4}\left[ \left(
b_{1q}+b_{2q}+b_{3q}+\beta \right) \exp \left( \mu _{q}t\right) \right] }~,
\end{equation}%
be \textit{time-independent}.

And since it is easily seen that
\end{subequations}
\begin{subequations}
\begin{equation}
f\left( t\right) =\frac{k_{I}\left( b_{31}+1\right) \varphi \left( t\right)
}{\left( b_{11}+b_{21}+b_{31}+\beta \right) }~,  \label{42ft}
\end{equation}%
with%
\begin{equation}
\varphi \left( t\right) =\frac{1+\sum_{q=2}^{4}\left\{ \left( \frac{b_{3q}+1%
}{b_{31}+1}\right) \exp \left[ \left( \mu _{q}-\mu _{1}\right) t\right]
\right\} }{1+\sum_{q=2}^{4}\left\{ \left( \frac{b_{1q}+b_{2q}+b_{3q}+\beta }{%
b_{11}+b_{21}+b_{31}+\beta }\right) \exp \left[ \left( \mu _{q}-\mu
_{1}\right) t\right] \right\} }~,  \label{42phit}
\end{equation}%
the requirement that $f\left( t\right) $ be \textit{time-independent}
amounts to the following $3$ constraints:
\end{subequations}
\begin{equation}
\frac{b_{3q}+1}{b_{31}+1}=\frac{b_{1q}+b_{2q}+b_{3q}+\beta }{%
b_{11}+b_{21}+b_{31}+\beta }~,~~~q=2,3,4~,
\end{equation}%
which clearly entail---via the expressions (\ref{41bel})---$3$ corresponding
\textit{constraints} on the parameters of the original model (\ref{eqs1}).

\bigskip

\section{Concluding remarks}

In this paper we have identified certain solutions of the \textit{pandemic}
model introduced in the paper \cite{RWK2021}; these solutions---and the
constraints on the parameters of the model required for their validity---are
all identified by \textit{algebraic} equations which can be \textit{%
explicitly} solved; we did not report the corresponding explicit formulas
because they are so complicated to be hardly useful when written for \textit{%
a priori arbitrary} assignments of the parameters of the pandemic model,
while they can instead be easily managed for any \textit{specific numerical}
assignment of these parameters. We therefore leave the utilization of these findings to
the interested pandemics experts.

Additional solutions---more special but perhaps displaying more interesting
evolutions---correspond to the special cases in which the algebraic
quartic-equation (\ref{eq8c}) features $4$ roots $\mu _{n}$ which are
\textit{not} all different among themselves. This case shall perhaps be
eventually treated in a separate paper by ourselves or others.

\bigskip

\section*{Acknowledgments}

We like to acknowledge with thanks $3$ grants, facilitating our collaboration---mainly
developed via e-mail exchanges---by making it possible for FP to visit three times
the Department of Physics of the University of Rome "La Sapienza": two
granted by that University, and one granted jointly by the Istituto
Nazionale di Alta Matematica (INdAM) of that University and by the
International Institute of Theoretical Physics (ICTP) in Trieste in the
framework of the ICTP-INdAM "Research in Pairs" Programme.\textbf{\ }%
Finally, we also\ like to thank Fernanda Lupinacci who, in these difficult
times---with extreme efficiency and kindness---facilitated all the
arrangements necessary for the presence of FP with her family in Rome.

\end{document}